\documentclass[runningheads]{llncs}
\usepackage{latexsym,bm}
\def\ra{\rightarrow}
\def\ral{\longrightarrow}
\def\be{\begin{equation}}
\def\ee{\end{equation}}
\def\bea{\begin{eqnarray}}
\def\eea{\end{eqnarray}}

\def\HH{{\cal H}}
\def\O{{\cal O}}
\def\F{{\bf F}}
\def\Z{{\bf Z}}
\def\C{{\bf C}}
\def\Q{{\bf Q}}
\def\R{{\bf R}}
\def\Pr{{\bf P}}
\def\Aut{{\rm Aut}}
\def\End{{\rm End}}
\def\SL{{\rm SL}}
\def\PSL{{\rm PSL}}
\def\PGL{{\rm PGL}}
\def\Gal{{\mathop{\rm Gal}}}
\def\X{{\rm X}}
\def\XX{{\cal X}}

\def\Xone{{\XX}(1)}

\def\XXoh{{\XX_0}}
\def\XXun{{\XX_1}}
\def\JJoh{{\cal J}_0}
\def\JJun{{\cal J}_1}
\def\0{^{\phantom0}}
\def\9{_{\phantom9}}
\def\wN{{w\0_N}}
\def\ww{{w\0_3}}
\def\Ovary{{\hbox{$\cal O$-variety}}}
\def\Ovaries{{\hbox{$\cal O$-varieties}}}

\def\Qed{~~$\diamondsuit$}
\begin{document}

\title{Shimura curves for level-3 subgroups\\
of the $(2,3,7)$ triangle group\\
\large\em and some other examples}
\titlerunning{Shimura curves $\XXoh(3)$, $\JJun(3)$
for the $(2,3,7)$ triangle group}
\author{Noam D.~Elkies}
\institute{Department of Mathematics, Harvard University, Cambridge, MA 02138 USA}
\maketitle

\begin{abstract}
The $(2,3,7)$ triangle group is known to be associated with
a quaternion algebra $A/K$\/ ramified at two of the three real places
of $K=\Q(\cos 2\pi/7)$ and unramified at all other places of~$K$.
This triangle group and its congruence subgroups thus give rise
to various Shimura curves and maps between them.
We study the \hbox{genus-$1$} curves $\XXoh(3)$, $\XXun(3)$
associated with the congruence subgroups $\Gamma_0(3)$, $\Gamma_1(3)$.
Since the rational prime~$3$ is inert in~$K$,
the covering $\XXoh(3)/\Xone$ has degree~$28$,
and its Galois closure $\XX(3)/\Xone$
has geometric Galois group $\PSL_2(\F_{27})$.
Since $\Xone$ is rational, the covering $\XXoh(3)/\Xone$
amounts to a rational map of degree~$28$.
We compute this rational map explicitly.
We find that $\XXoh(3)$ is an elliptic curve
of conductor~$147=3\cdot 7^2$ over~$\Q$,
as is the Jacobian $\JJun(3)$ of $\XXun(3)$;
that these curves are related by an isogeny of degree $13$;
and that the kernel of the \hbox{$13$-isogeny} from~$\JJun(3)$
to $\XXoh(3)$ consists of \hbox{$K$-rational points.}
We also use the map $\XXoh(3) \ra \Xone$ to locate
some complex multiplication (CM) points on~$\Xone$.
We conclude by describing analogous behavior of a few Shimura curves
associated with quaternion algebras over other cyclic cubic fields.
\end{abstract}
\section{Introduction}

\subsection{Review: Quaternion algebras over~$K$,
Shimura modular curves, and the $(2,3,7)$ triangle group.}

Let $K$\/ be the field $\Q(\cos 2\pi/7)$, which is the totally real
cubic field of minimal discriminant (namely, discriminant~$49$);
let $O_K$ be the ring $\Z[2\cos 2\pi/7]$ of algebraic integers in~$K$\/;
and let $A$ be a quaternion algebra over~$K$\/
ramified at two of the three real places and at no other place of~$K$.
This algebra exists because the set of ramified places has even
cardinality, and is determined uniquely up to $K$-algebra isomorphism.
(See for instance \cite{V} for this and other basic results
on quaternion algebras and Shimura curves,
and \cite{Shimcomp} for our computational context.)
All maximal orders in~$A$\/ are conjugate because $A$\/ is indefinite
and $O_K$ has narrow class number~$1$;
we fix one maximal order $\O \subset A$.
Let $\Gamma(1)$ be the group of elements of~$\O$ of reduced norm~$1$.
Since $A$ has exactly one unramified real place, $\Gamma(1)$
embeds into $\SL_2(\R)$ as a discrete \hbox{co-compact} subgroup.
Let $\HH$ be the upper half-plane, with the usual action of
$\PSL_2(\R) = \SL_2(\R) / \{\pm1\}$.  The quotient of $\HH$ by
$\Gamma(1)/\{\pm1\}$ then has the structure of a compact Riemann surface,
and thus of an algebraic curve over~$\C$.  In fact this quotient
has a natural structure as a curve over~$K$,
namely the Shimura curve associated to $\Gamma(1)$.
We call this Shimura curve $\Xone$, in analogy with the modular elliptic curve
$\X(1)$ of the classical theory of elliptic and modular curves:
as $\X(1)$ parametrizes elliptic curves, $\Xone$ parametrizes certain
abelian varieties which we shall call ``\Ovaries''.\footnote{
\label{footnote:carayol}
  Warning: unlike the case of Shimura curves associated with
  quaternion algebras over~$\Q$, here an \Ovary\ is not simply
  a principally polarized abelian variety with endomorphisms by~$\O$.
  Indeed there can be no abelian variety~$V$\/
  with $\End(V) \otimes \Q = A$,
  because the set of ramified primes of~$A$\/
  neither contains nor is disjoint from the set of real places of~$K$.
  (See for instance \cite{Mumford}, specifically Thm.1 on p.192
  (positivity of the Rosati involution) together with
  the classification in Thm.2 on p.201.)
  The moduli description is quite complicated
  (see \cite{Shimura,Carayol}),
  and requires an auxiliary quadratic extension $K'/K$,
  with the field $K'$ totally imaginary; for instance,
  for our case $K=\Q(\cos 2\pi/7)$ we may take
  $K'=K(\sqrt{-7})=\Q(e^{2\pi i/7})$.
  The moduli description yields $\Xone$ as a curve over~$K'$,
  but fortunately $\Xone$ descends to a curve over~$K$\/
  independent of the choice of~$K'$.
  The same is true of the curves $\XX(N), \XXun(N), \XXoh(N)$
  and the maps $\XX(N) \ra \XXun(N) \ra \XXoh(N) \ra \Xone$
  (to be introduced soon) that concern us in this paper,
  at least provided that the extension $K'/K$\/
  satisfies some local conditions at the primes dividing~$N$.
%
%
  I am grateful to Benedict Gross for explaining these subtleties
  and suggesting the references to~\cite{Carayol,Mumford}.
  Fortunately very little of this difficult theory
  is needed for our computations.
  }
By work of Shimura~\cite{Shimura}, based on the classical work
of Fricke \cite{F1,F2}, the group
$\Gamma(1)/\{\pm1\} \subset \SL_2(\R)$ is the $(2,3,7)$ triangle group
(the group generated by products of pairs of reflections in
the sides of a hyperbolic triangle of angles $\pi/2$, $\pi/3$, $\pi/7$).
Hence $\Xone$ is a curve of genus~$0$
with elliptic points of orders $2$, $3$,~$7$.
We choose the rational coordinate
$t : \Xone \stackrel{\sim}{\ral} \Pr^1$
that takes the values $0,1,\infty$ respectively at these three points;
this determines $t$ uniquely (as the classical modular invariant $j$
is determined by its values $1728,0,\infty$ at the elliptic points
of orders $2$, $3$ and the cusp).

Suppose now that $\Gamma \subset \Gamma(1)/\{\pm1\}$
is a subgroup of finite index~$d$.
Then $\XX = \HH/\Gamma$ is a compact Riemann surface
with a \hbox{degree-$d$}\/ map to~$\Xone$
that is unramified away from the elliptic points of~$\Xone$.
Composing the map $\pi: \XX \ra \Xone$ with our isomorphism
$t : \Xone \stackrel{\sim}{\ral} \Pr^1$
yields a rational function $t \circ \pi$ of degree~$d$\/ on~$\XX$\/
that is unramified except above $t=0$, $t=1$, and $t=\infty$.
Such a rational function is often called a ``Belyi function'' on~$\XX$,
in tribute to Belyi's theorem~\cite{Belyi} that a Riemann surface
admits such a function if and only if it can be realized as
an algebraic curve over $\overline\Q$.  It will be convenient,
and should cause no confusion,
to use $t$ also as the notation for this rational function on~$\XX$.

In particular, if $\Gamma$ is a congruence subgroup of $\Gamma(1)/\{\pm1\}$
--- that is, if there exists a nonzero ideal~$N$\/ of~$O_K$
such that $\Gamma$ contains the image in $\Gamma(1)/\{\pm1\}$
of the normal subgroup
$$
\Gamma(N) := \{ g \in \Gamma(1)  \mid  g \equiv 1 \bmod N\O \}
$$
--- then $\XX$\/ is also a Shimura modular curve,
parametrizing \Ovaries\ with some \hbox{level-$N$}\/ structure.
For example, $\Gamma(N)$ itself is a congruence subgroup,
corresponding to a Shimura curve
parametrizing \Ovaries\ with full \hbox{level-$N$}\/ structure.
We call this curve $\XX(N)$, again in analogy with the classical case.
If $N$\/ is a prime ideal with residue field~$k$,
the ring $\O/N\O$ is isomorphic with the ring $M_2(k)$
of $2 \times 2$ matrices over~$k$,
because $A$\/ is unramified at all finite places of~$k$, including~$N$.
We choose one isomorphism $\iota: \O/N\O \stackrel{\sim}{\ral} M_2(k)$;
the choice does not matter
because all are equivalent under conjugation by~$\O^*$.
%
%
We then obtain congruence subgroups
$\Gamma_0(N) \rhd \Gamma_1(N) \supset \Gamma(N)$,
where $\Gamma_0(N)$ consists of those $g \in \Gamma(1)$
for which $\iota(g)$ is upper triangular, and $\Gamma_1(N)$
is the kernel of the map from~$\Gamma_0(N)$ to $k^*/\{\pm1\}$
taking $g$ to the top left entry of the matrix $\iota(g)\bmod N\O$.
The corresponding Shimura curves are naturally denoted
$\XXoh(N)$, $\XXun(N)$ respectively.
Once one has formulas for the cover $\XXoh(N) \ra \Xone$,
one can use them as in~\cite{Serre}
to obtain explicit extensions of~$K$\/ or~$\Q$
(the latter when the residual characteristic of~$N$\/
does not split in~$O_K$) whose normal closures have Galois groups
$\PSL_2(k)$, $\Aut(\PGL_2(k))$,
or other groups intermediate between these two.

Analogous to the Atkin-Lehner involution $\wN$ of the classical modular
curve $\X_0(N)$, we have an involution $\wN$ of $\XXoh(N)$,
which can be constructed either from the normalizer of $\Gamma_0(N)$ in~$A^*$
or by invoking the dual isogeny of the ``cyclic'' isogeny between \Ovaries\
parametrized by a generic point of~$\Gamma_0(N)$.
Once explicit formulas are known for both the map $\XXoh(N) \ra \XX(1)$
and the involution $\wN$, one can easily compute other interesting data.
By eliminating $p\in\XXoh(N)$ from the system
$T_1=t(p)$, $T_2=t(\wN(p))$, we obtain
the ``modular polynomial'' $\Phi_N(T_1,T_2)$,
which is the irreducible polynomial such that $\Phi_N(t_1,t_2)=0$
if and only if $t_1,t_2\in\C$ are \hbox{$t$-coordinates}
of \hbox{``$N$-isogenous''} points on~$\Xone$.
The solutions of $\Phi_N(t,t)=0$
are then coordinates of points of complex multiplication (CM)
on~$\Xone$, that is, points that parametrize \Ovaries\ with~CM.
Finally, with some more effort we can obtain equations defining
the recursive tower of curves $\XXoh(N^n\9)$ ($n=1,2,3,\ldots$),
whose reduction modulo any prime~$\varpi$ of~$O_K$ not dividing~$N$\/
%
%
%
%
is asymptotically optimal over the quadratic extension
of the residue field of~$\varpi$ (see~\cite{Modtower}).

If $\Gamma$ is a proper normal subgroup, congruence or not,
in $\Gamma(1)/\{\pm1\}$, then the quotient group~$G$\/ acts on~$\XX$.
By the Riemann-Hurwitz formula, $d=|G|=84(g-1)$,
where $g$ is the genus of~$\XX$.
That is, $(\XX,G)$ is a ``Hurwitz curve'':
a curve of genus $g>1$ that attains the Hurwitz bound
(\cite{Hurwitz}, see also \cite[Ch.I, Ex.F-3 ff., pp.45--47]{ACGH})
on the number of automorphisms
of a curve of genus~$g$ over a field of characteristic~zero.
Conversely, all Hurwitz curves arise in this way
for some proper normal subgroup of~$\Gamma(1)$.
For example, if $N$\/ is a prime ideal with residue field~$k$\/ then
$\XX(N)$ is a curve with $G=\PSL_2(k)$ that attains the Hurwitz bound.

The first few possibilities for~$N$\/ yield
$k$ of size $7$, $8$, $13$, and~$27$.  The first of these, with
$N$\/ the prime of~$K$\/ above the totally ramified rational prime~$7$,
yields the famed Klein quartic of genus~$3$,
which also arises as the classical modular curve $\X(7)$
(see \cite{NDE:Klein}).
The second, with $N$\/ above the inert rational prime~$2$,
yields the Fricke-Macbeath curve of genus~$7$ \cite{FrickeM,Mac}.
If $N$\/ is any of the three primes above the split rational prime~$13$
then $\XX(N)$ is a curve of genus~$14$ that can be defined over~$K$\/
but not over~$\Q$~\cite{Streit}.
In each of those cases, $\XXoh(N)$ has genus~zero.

The next case, and the main focus of this paper,
is the prime of residue field $\F_{27}$ above the inert rational prime~$3$.
We call the resulting Shimura curves $\XXoh(3)$, $\XXun(3)$, $\XX(3)$.
%
%
%
%
%
%
Since this prime is invariant under $\Gal(K/\Q)$,
these Shimura curves and their Belyi maps to~$\Xone$
(though not the action of $\PSL_2(\F_{27})$ on $\XX(3)$)
can be defined over~$\Q$~\cite{DN,Wolfart}.
Here for the first time $\XXoh(N)$ has positive genus;
we calculate that $\XXoh(3)$ and $\XXun(3)$ both have genus~$1$.
The determination of the curve $\XXoh(3)$,
and of its \hbox{degree-$28$} Belyi map to $\XX(1)$ ---
that is, the \hbox{degree-$28$} function $t$\/ on~$\XXoh(3)$ ---
requires techniques beyond those of~\cite{Shimcomp}.

In section~2 of this paper,
we exhibit equations for this cover and show how we compute them.
We locate the coefficients to high \hbox{$p$-adic} precision
via Newton's method, using $p=29$, a prime at which the cover
has good reduction but some of the coefficients are known in advance
modulo~$p$ and the rest can be determined algebraically.
We then compute the simplest rational numbers consistent
with the \hbox{$p$-adic} approximations,
and prove that they are in fact the correct coefficients.
The same method, possibly extended by exhaustive searching mod~$p$,
can be used to compute other such covers.  (For instance,
we used it to compute a previously unknown Shimura-Belyi function
of degree~26 connected with the $(2,3,8)$ triangle group.
We shall give the details of this computation,
and other Shimura-Belyi maps for the $(2,3,8)$ triangle group,
in a later paper.)
%
%
%
%
%
%
%
%
%
%

Having obtained equations for the map $\XXoh(3) \ra \XX(1)$, we use
those equations as explained above to locate several CM~points on~$\Xone$.
In particular, we confirm our conjecture from~\cite{Shimcomp}
for the rational point of CM~discriminant~$-11$: we had computed
its \hbox{$t$-coordinate} to high enough precision to convincingly
guess it as an element of~$\Q$, but could not prove this guess.
It also follows immediately from our formulas that
$\XXoh(3)$ has a \hbox{$\Q$-rational} point, and is thus an elliptic curve.
In section 3 of this paper, we note some arithmetic properties
of this curve, and use the modular structure to explain them.
In particular, $\XXoh(3)$ admits a rational \hbox{$13$-isogeny}
to another elliptic curve,
which we identify as the Jacobian of $\XXun(3)$.
We then also explain why the kernel of the dual isogeny from this
Jacobian to $\XXoh(3)$ must consist of \hbox{$K$-rational} points,
a fact first noted via direct computation by Mark Watkins~\cite{Watkins}.

In the final section we give some other examples
suggested by these ideas.
Watkins later found in the tables of elliptic curves
two pairs of \hbox{$7$-isogenous} curves,
of conductors $162=2\cdot 9^2$ and $338=2\cdot 13^2$,
with \hbox{$7$-torsion} structure over the cyclic cubic fields
of discriminant $9^2$ and $13^2$.
%
%
%
%
%
He suggested that each of these pairs might thus be the Jacobians
of the Shimura modular curves $\XXoh(2), \XXun(2)$
for a quaternion algebra over the cyclic cubic field
of discriminant $9^2$, $13^2$ respectively
that is ramified at two of its three real places
and at no finite primes.
We verify that this suggestion is correct
for the curves of conductor~$2\cdot 9^2$,
corresponding to a subgroup of the $(2,3,9)$ triangle group.
We cannot at present compute the analogous curves and maps
for the Shimura curves of level~$2$ associated with a quaternion algebra
over the cubic field of conductor~$13^2$,
though we verify that both curves have genus~$1$.
It might be possible to identify the curves using more refined
arithmetical information as in \cite{GonRot,Kurihara},
but this would still leave open the problems of explicitly computing
the maps $\XXoh(2)\ra\Xone$ and $w_2$.

\section{The curve $\XXoh(3)$}

\subsection{Preliminaries}

We use the notations of the previous section.
Let $t$ be the \hbox{degree-$28$} Belyi function on $\XXoh(3)$.
The elements of orders $2$, $3$, $7$ in $\PSL_2(\F_{27})$
act on $\Pr^1(\F_{27})$ by permutations with cycle structure
$2^{14}$, $3^9 1$, $7^4$.
We denote by $P_0$ the unique simple zero of $(t-1)$ on $\XXoh(3)$
(corresponding to the unique fixed point of an element of order~$3$
in $\PSL_2(\F_{27})$).
Since $\XXoh(3)$ and $t$ are defined over~$\Q$, so is~$P_0$.
The preimages other than~$P_0$ of $0,1,\infty$ under~$t$
constitute disjoint effective divisors $D_{14}$, $D_9$, $D_4$, which
are the sum of $14$, $9$, $4$ distinct points on $\XXoh(3)$ respectively,
such that
$$
t^*(0) = 2 D_{14},
\quad
t^*(1) = 3 D_9 + (P_0),
\quad
t^*(\infty) = 7 D_4.
$$
By the Riemann-Hurwitz formula, $\XXoh(3)$ is a curve of genus~$1$.
We give $\XXoh(3)$ the structure of an elliptic curve
by choosing $P_0$ as the origin of the group law.
Since $P_0$ is the unique point of~$\XXoh(3)$
parametrizing a cyclic \hbox{$3$-isogeny}
between the \hbox{order-$3$} elliptic point $t=1$ of~$\Xone$ and itself,
$P_0$ must be fixed under~$\ww$.  Therefore $\ww$ is the unique
involution of $\XXoh(3)$ as an elliptic curve, namely,
multiplication by~$-1$.

{\bf Proposition 1.}
{\em
i) The differential $dt$\/ has divisor $D_{14} + 2 D_9 - 8 D_4$.
The divisors $D_{14}$, $D_9$, $D_4$
are linearly equivalent to $14(P_0)$, $9(P_0)$, $4(P_0)$, respectively.
\\
$\phantom{|}$\quad ii)
There are nonzero rational functions $F_{14}$, $F_9$, $F_4$ on~$\XXoh(3)$
with divisors $D_{14}-14(P_0)$, $D_9-9(P_0)$, $D_4-4(P_0)$.
For each choice of $F_{14}, F_9, F_4$,
there exist nonzero scalars $\alpha,\beta$ such that
$F_{14}^2 = \alpha F_4^7 + \beta F_9^3$ and $t=F_{14}^2/\alpha F_4^7$.
If the $F_n$ are defined over~$\Q$
then $\alpha,\beta$ are rational as well.
}

{\em Proof}\/: i) Since $\XXoh(3)$ has genus~one, the divisor of~$dt$\/
is linearly equivalent to zero.  This divisor is regular
except for poles of order~$8$ at the four points of~$D_4$.
Moreover, $dt$\/ has simple zeros at the points of $D_{14}$
and double zeros at the points of~$D_9$.
This accounts for $14 + 2\cdot 9 = 32$ zeros,
same as the number $8\cdot 4$ of poles, so $dt$\/ has no further zeros
(which we could also have deduced from the fact that the map~$t$\/
is unramified when $t\notin\{0,1,\infty\}$).

It follows that $D_{14} + 2 D_9 \sim 8 D_4$.
This together with the linear equivalence of the divisors
$t^*(0)$, $t^*(1)$, $t^*(\infty)$ yields our claim that
$D_n \sim n (P_0)$ for each $n=14,9,4$; for instance
\begin{eqnarray*}
D_4 & = & 49 D_4 - 48 D_4 \sim 7 (3D_9+(P_0)) - 6(D_{14}+2D_9)
\\
&& = 3 (3D_9-2D_{14}) + 7 (P_0) \sim -3 (P_0) + 7 (P_0) = 4(P_0).
\end{eqnarray*}
ii) Functions $F_n$ with divisors $D_n - n(P_0)$ exist by part~(i).
The rational functions $t$\/ and $F_{14}^2/F_4^7$ have the same divisor,
so $F_{14}^2/F_4^7 = \alpha t$ for some nonzero scalar~$\alpha$.
Likewise $t-1$ and $F_9^3/F_4^7$ have the same divisor, so
$F_9^3/F_4^7 = \alpha_0(t-1)$ for some nonzero scalar~$\alpha_0$.
Eliminating $t$ yields the desired identity
$F_{14}^2 = \alpha F_4^7 + \beta F_9^3$ with $\beta = \alpha/\alpha_0$.
It follows that $\alpha$, and thus also $\beta$, is rational
if each $F_n$ is defined over~$\Q$.
This completes the proof of Proposition~1.\Qed

Now fix a nonzero holomorphic differential $\omega$ on~$\XXoh(3)$
defined over~$\Q$, and define a derivation $f\mapsto f'$
on the field of rational functions on~$\XXoh(3)$ by $df = f' \omega$.

{\bf Proposition 2.}
{\em
i) Let $F_{14}, F_9, F_4$ be as in Proposition~1(ii).
Then $F_{14}\0$, $F_9^2$, $F_4^6$ are scalar multiples of
$$
3 F_4\0 F_9' - 7 F_4' F_9\0,
\quad
2 F_4\0 F_{14}' - 7 F_4' F_{14}\0,
\quad
2 F_9\0 F_{14}' - 3 F_9' F_{14}\0
$$
respectively.
\\
ii) The functions
$$
\frac{14 F_9\0 F_4'' - 13 F_9' F_4'}{F_4},
\quad
\frac{6 F_4\0 F_9'' - 29 F_4' F_9'}{F_9},
\quad
\frac{6 F_4\0 F_{14}'' - 43 F_4' F_{14}'}{F_{14}}
$$
on~$\XXoh(3)$ are regular except at $P_0$.
}

{\em Proof}\/:
i) By Proposition~1 we have $t = F_{14}^2 / \alpha F_4^7$.
Taking the logarithmic derivative, we find
$$
\frac{dt}{t}
= \left( 2 \frac{F_{14}'}{F_{14}} - 7 \frac{F_4'}{F_4} \right) \omega
= \frac{2 F_4\0 F_{14}' - 7 F_4' F_{14}\0}{F_4 F_{14}} \, \omega.
$$
Since $dt$\/ has divisor $D_{14} + 2 D_9 - 8 D_4$,
the divisor of $dt/t$\/ is $2 D_9 - D_{14} - D_4$.
Thus $F_4 F_{14} dt/t$ is a differential with divisor
$2 D_9 - 18 (P_0)$, same as the divisor of $F_9^2$.
Since the divisor of a differential on a genus-one curve
is linearly equivalent to zero, it follows that
$2 F_4\0 F_{14}' - 7 F_4' F_{14}\0$ is a multiple of $F_9^2$,
as claimed.
The formulas for $F_{14}$ and $F_4^6$ are obtained in the same way
by computing the logarithmic derivatives of $t-1$ and $t/(t-1)$
respectively.

ii) Each of these is obtained by substituting
one of the identities in~(i) into another.
Since the computations are similar and we shall use
only the result concerning $(6 F_4\0 F_9'' - 29 F_4' F_9')/F_9\0$,
we prove this result, and again leave the other two as exercises.

We have $F_{14} = \gamma_1 (3 F_4\0 F_9' - 7 F_4' F_9\0)$
and $F_9^2 = \gamma_2(2 F_4\0 F_{14}' - 7 F_4' F_{14}\0)$
for some nonzero scalars~$\gamma_1, \gamma_2$.  
Hence
\begin{eqnarray*}
&
F_9^2 = \gamma_1 \gamma_2 \bigl( 2 F_4\0 (3 F_4\0 F_9' - 7 F_4' F_9\0)'
- 7 F_4' (3 F_4\0 F_9' - 7 F_4' F_9\0) \bigr)
&
\\
&
\qquad = \gamma_1 \gamma_2 \bigl( (6 F_4\0 F_9'' - 29 F_4' F_9') F_4\0
+ (49 {F_4'}^2 - 14 F_4\0 F_4'') F_9\0 \bigr) \, .
&
\end{eqnarray*}
Dividing by $F_9$, we find that
$$
(6 F_4\0 F_9'' - 29 F_4' F_9') \frac{F_4}{F_9}
+ (49 {F_4'}^2 - 14 F_4\0 F_4'') = \frac1{\gamma_1 \gamma_2} F_9,
$$
which is regular away from~$P_0$.
The same is true for $49 {F_4'}^2 - 14 F_4\0 F_4''$,
and thus also for $(6 F_4\0 F_9'' - 29 F_4' F_9') F_4\0 / F_9\0$.
Since $F_4$ and $F_9$ have no common zeros, it follows that
$(6 F_4\0 F_9'' - 29 F_4' F_9')/F_9$ has no poles except for
a multiple pole at $P_0$, as claimed.\Qed

\subsection{Computation}

{\bf Theorem 1.}
{\em
The curve $\XXoh(3)$ is isomorphic with the elliptic curve
$$
y^2 + y = x^3 + x^2 - 44704 x - 3655907.
$$
The functions $F_n$ on this curve may be given by
\begin{eqnarray*}
F_4
& = &
  x^2 - 1208 x - 227671 + 91 y,
\\
F_9
& = &
(8 x^3 - 384 x^2 - 13656232 x - 678909344) y
\\
&& {} \!\!\! -
(1015 x^4 + 770584 x^3 + 163098512 x^2 + 29785004488 x + 2319185361392),
\\
F_{14}
& = &
8 x^7 + 400071 x^6 - 343453068 x^5 - 238003853192 x^4
\\
&&
{} - 116011622641292 x^3 - 15704111899877744 x^2
\\
&&
{} - 95316727595264672 x + 53553894620234333456
\\
&& {} \!\!\! -
(8428 x^5 + 19974360 x^4 + 18880768004 x^3 + 4128079708928 x^2
\\
&&\quad
{} + 335969653582304 x + 17681731246686360) y,
\end{eqnarray*}
satisfying
$$
64 F_4^7 - 343 F_9^3 = F_{14}^2,
$$
and then
$$
t = \frac{F_{14}^2}{64 F_4^7} = 1 - \frac{343 F_9^3}{64 F_4^7} \, .
$$
}

{\em Proof}\/:
With a symbolic algebra package one may readily confirm
the identity among the $F_n$, and might even feasibly
(albeit not happily) verify the Galois group by following
the $28$ preimages of a point on the \hbox{$t$-line}
as that point loops around $t=0$ and $t=1$;
this would suffice to prove the theorem
(since the cover $\XXoh(3) \ra \XX(1)$ is determined
by its Galois group and ramification behavior),
but would not explain the provenance of the formulas.
We thus devote most of the proof to the computation of the~$F_n$.

We begin by observing that our proofs of Propositions~1 and~2
used the ramification behavior of the cover $\XXoh(3) \ra \XX(1)$,
but not its Galois group or the Shimura-curve structure.
This will remain true in the rest of our computation.\footnote{
  We shall retain the notation $\XXoh(3)$ for the curve
  and $\ww$ for the involution, rather than introduce
  separate notations for an elliptic curve
  that we have not yet identified with $\XXoh(3)$
  and the multiplication-by-$(-1)$ map on the curve.
  }
We thus show that the ramification behavior uniquely determines
the \hbox{degree-$28$} cover.  In particular this yields
the following purely group-theoretical statement
(that we could also have checked directly):
any permutations $\sigma_2, \sigma_3, \sigma_7$ of $28$ letters
with cycle structure $2^{14}$, $3^9 1$, $7^4$ such that
$\sigma_2 \sigma_3 \sigma_7 = 1$ are equivalent under conjugation
in $S_{28}$ with our generators of $\PSL_2(\F_{27})$,
and that we do not have to verify the Galois group
as suggested in the previous paragraph.

Our strategy is to first find the cover modulo the prime~$29$,
which occurs in one of the formulas in Proposition~2(ii),
and then use Newton's method over $\Q_{29}$
to compute a lift of the coefficients
to sufficient \hbox{$29$-adic} precision
to recognize them as rational numbers.
The prime $29$ is large enough to guarantee good reduction
of any Belyi cover of degree~$28$:
if a Belyi cover has bad reduction at some prime~$p$
then the normal closure of the cover has a Galois group whose order
is a multiple of~$p$~\cite{Beckmann}.

We may assume that each $F_n$ ($n=14,9,4$) is scaled so that
it has \hbox{$29$-adically} integral coefficients
and does not reduce to zero mod~$29$.
In characteristic $29$, the second function in Proposition~2(ii)
simplifies to $6 F_4\0 F_9''/F_9\0$.  Again we use the fact that
$F_4,F_9$ have no common zeros to conclude that
$\xi := F_9''/F_9\0$ is regular away from~$P_0$.
At $P_0$, we know that $F_9$ has a pole of order~$9$;
thus $F_9''$ has a pole of order~$11$, and $\xi$ has a double pole.
Since $\xi$ is regular elsewhere,
it follows that $\xi$ is invariant under~$\ww$.

We claim that $F_9$ is anti-invariant, that is, that $w_3(F_9)=-F_9$.
Indeed the differential equation $f''=\xi f$ satisfied by~$F_9$
must also hold for the invariant and anti-invariant parts of~$F_9$,
call them $F_9^+, F_9^-$.
Now $F_9^-$ has a pole of order~$9$ at~$P_0$.
If $F_9^+$ is not identically zero then its valuation at~$P_0$
is $-d$\/ for some $d\in\{0,2,4,6,8\}$;
comparing leading coefficients in the local expansion about~$P_0$ of
${F_9^+}{}'' = \xi F_9^+$ and ${F_9^-}{}'' = \xi F_9^-$,
we obtain $9\cdot 10 \equiv d(d+1) \bmod 29$, which is impossible.
Hence $F_9^+ = 0$ as claimed.
(We could also have obtained $F_9^+ = 0$ by arguing that $t=1$
is the only supersingular point on $\Xone$ mod~$29$,
whence the zeros of $F_9$, which are the preimages of $t=1$,
must be permuted by~$\ww$.  But this would break our promise
not to rely on the modular provenance of the cover.)

Now let $y^2=x^3+ax+b$ be a (narrow) Weierstrass equation
for $\XXoh(3)$, and choose $\omega=dx/y$.
Then our derivation $f \mapsto f'$ is characterized by
$x'=y$ and $y'=(3x^2+a)/2$.  For the rest of this paragraph
we calculate modulo~$29$.  We have seen that $F_9$
is a scalar multiple of $(x^3 + s_1 x^2 + s_2 x + s_3) y$
for some $s_1,s_2,s_3 \in \overline{\F_{\!29}}$.  Equating coefficients
in $F_9'' = \xi F_9\0$, we find first that $\xi = 8 x + 6 s_1$,
then that $s_2\0 = -12 s_1^2 - 8a$ and
$s_3\0 = 7 b - 3 s_1^3 - a s_1\0$,
and finally that $s_1 b + s_1^4 + 9 a s_1^2 + 9 a^2 = 0$.
We have also seen that $F_4$ is a scalar multiple of
$x^2 + t_1 y + t_2 x + t_4$
for some $t_1,t_2,t_4 \in \overline{\F_{\!29}}$.
Using Proposition~2(i) we compute $F_{14}$ and $F_4^6$ up to scaling
in terms of $s_1,a,b$ and the $t_i$, and compare with
$(x^2 + t_1 y + t_2 x + t_4)^6$.
After matching the leading (\hbox{degree-$24$}) coefficients,\footnote{
  The ``degree'' of a polynomial in~$x,y$
  is the order of its pole at~$P_0$.  The vector space
  of rational functions on~$\XXoh(3)$ that are regular away from~$P_0$
  has basis $\{x^i|i=0,1,2,\ldots\} \cup \{x^i y|i=0,1,2,\ldots\}$;
  the monomials $x^i$, $x^i y$ have degrees $2i$, $2i+3$.
  }
we find that the \hbox{degree-$23$} coefficients agree identically,
but in degree~$22$ we find $t_2 = 11 (t_1^2-s_1\0)$,
and the degree~$21$ comparison yields $t_1=0$ or $t_1^2 = 5 s_1\0$.
We cannot have $t_1=0$,
for then $F_4$ would be even, and then so would $F_{14}$
(since $F_{14}$ is proportional to $3 F_4\0 F_9' - 7 F_4' F_9\0$),
which is impossible since the nonzero odd function $F_9^3$
is a linear combination of $F_4^7$ and $F_{14}^2$.
Therefore $t_1^2 = 5s_1$.
Comparing the next few coefficients in our two expressions for $F_4^6$,
we learn that $t_4 = s_1^2 + 3 a$ and $a = 9 s_1^2$.
This completes the determination of our cover mod~$29$ up to scaling.
For instance, we may use the equation
$y^2 = x^3 + 9x + 1$ for $\XXoh(3)$ over $\overline{\F_{\!29}}$,
and then check that
$$
F_4 = x^2 + 11 y - 14 x - 1,
\quad
F_9 = (x^3 + x^2 + 3 x - 5) y,
$$ \vspace*{-3ex} $$
F_{14} = (x^7-14x^6-5x^5-9x^4-10x^3+2x^2+10x-7)
- (8x^5-x^3-3x^2+x+3)y
$$
satisfy the identity $F_4^7 - 6 F_9^3 = F_{14}^2$.

We now use Hensel's lemma to show
that there is a unique such identity over $\Z_{29}$.
To do this we regard $\alpha F_4^7 + \beta F_9^3 = F_{14}^2$
as an overdetermined system of equations
in $a,b,\alpha,\beta$, and the coefficients of $F_4$ and $F_9$.
More precisely, we eliminate the ambiguity in the various choices
of scalar multiples by requiring that $F_4$ and $F_9$
each have leading coefficient~$1$, setting $b=9a$,
and defining $F_{14}$ by the formula
$3 F_4\0 F_9' - 7 F_4' F_9\0$ from Proposition~2(i).
(Any elements of~$\Z_{29}$ that reduce to~$1$ mod~$29$
would do for the leading coefficients, as would
any element of~$\Z_{29}$ that reduces to~$9$ mod~$29$ for $b/a$.)
We include in our equations the $x^i$ coefficients of~$F_9$
($i=0,1,2,3,4$), which were known to vanish mod~$29$
but not in~$\Z_{29}$.  By our analysis thus far,
this system of equations has a unique solution mod~$29$.
We compute the Jacobian matrix of our system of equations
at that solution, and verify that this matrix has full rank.
Therefore the solution lifts uniquely to $\Z_{29}$ by Hensel's lemma.
This is also the unique solution over~$\Q_{29}$, because
$29$ is a prime of good reduction as noted above.

It remains to recover the coefficients as rational numbers.
We compute them mod~$29^{128}$ by applying a Newton iteration
$7$~times \hbox{$29$-adically}.
(At each step, instead of computing the derivative mod~$29^{2^{n-1}}$
symbolically, we approximate each partial derivative by evaluating
the function mod~$29^{2^n}$ at two points that differ by
$29^{2^{n-1}}$ times the corresponding unit vector.
We also streamline the computation by requiring at each step
that only the $x^i$ coefficients match, ignoring the $x^i y$ terms;
this suffices because the resulting submatrix of the Jacobian matrix
still has full rank mod~$29$.)
We then use \hbox{$2$-dimensional} lattice reduction
%
%
%
%
%
to guess the rational numbers represented by those
\hbox{$29$-adic} approximations, and verify the guess
by checking that the resulting $F_4^7$, $F_9^3$, $F_{14}^2$
are \hbox{$\Q$-linearly} dependent.  This completes the determination
of the cover $\XXoh(3)\ra\Xone$ and the proof that it is characterized
by its degree and ramification behavior.
Finally we bring the curve $y^2=x^3+ax+b$ to minimal form
%
%
%
%
and replace $F_4, F_9, F_{14}$ by the smallest rational multiples
that eliminate the denominators of those rational functions.
This yields the formulas in Theorem~1.\Qed

In the same way that we proved $F_9$ is an odd function mod~$29$,
we can use the other two formulas in Proposition 2(ii) to prove
that $F_4$ and $F_{14}$ are even mod~$13$ and~$43$ respectively.
This is confirmed by Theorem~1:
the terms containing $y$ in $F_4$ and $F_{14}$
are divisible by~$13$ and~$43$ respectively.
The corresponding result for $F_9 \bmod 29$
is obscured by the minimal form of $\XXoh(3)$,
which makes it harder to recognize odd functions;
but $F_9$ can be seen to be congruent mod~$2y+1$
to a multiple of~$29$, as expected.

\subsection{Application: the CM point of discriminant $\bm{-11}$}

In~\cite[5.3]{Shimcomp}, we noted that $\Xone$ has a unique
CM~point of discriminant $-11$, and therefore that this point
has rational {$t$-coordinate}.  We then described
high-precision numerical computations
that strongly suggest this {$t$-coordinate} is
$$
\frac{88983265401189332631297917}{45974167834557869095293} =
\frac{7^3 43^2 127^2 139^2 307^2 659^2 11} {3^3 13^7 83^7},
$$
with $t-1$ having numerator $2^9 29^3 41^3 167^3 281^3$.
But we could not prove that this is correct.
Our formulas for the cover $\XXoh(3)\ra\Xone$ now let us do this:

{\bf Corollary to Theorem 1.}
{\em
The\/ {\rm CM-11} point on $\Xone$ has $t$-coordinate equal
$7^3 43^2 127^2 139^2 307^2 659^2 11 / 3^3 13^7 83^7$,
and lies under the two points of $\XXoh(3)$ with
$(x,y) = (-10099/64, -1/2 \pm 109809 \sqrt{-11} / 512)$.
}

{\em Proof}\/: The CM-11 point has two preimages in $\XXoh(3)$
that are switched by~$\ww$, corresponding to a pair of
``\hbox{$3$-isogenies}'' between
(the abelian variety parametrized by) that point and itself,
namely the pair $(-1 \pm \sqrt{-11})/2$ of \hbox{norm-$3$} elements
of the endomorphism ring.
A point of $\Xone$ is \hbox{$3$-isogenous} to itself
if and only if its {$t$-coordinate} satisfies $t = t(P) = t(\ww P)$
for some point~$P$\/ on~$\XXoh(3)$.  Equivalently, $P$\/ is
either a common pole of the functions $t(P), t(\ww P)$ on~$\XXoh(3)$
or a zero of the rational function $t(P) - t(\ww P)$ on that curve.
The former cannot happen here because by our formulas
$F_4$ and $F_4 \circ \ww$ have no common zeros.
Since $t(P) - t(\ww P)$ is an odd function on~$\XXoh(3)$,
we easily locate its zeros by writing it as $2y+1$
times a rational function of~$x$.
We find that the latter function vanishes at
$x = \infty$, $x = -1097/8$, $x = -10099/64$,
and the roots of four irreducible polynomials of degree~$3$
and two irreducibles of degree~$6$.
At $x = \infty$ we have $P=P_0$ and $t=1$.
At $x = -1097/8$ we have $y = -1/2 \pm 6615 \sqrt{-2} / 32$,
and calculate that $t = 1092830632334/1694209959$, which we already
identified as the \hbox{CM-8} point by a similar computation
with the curve $\XXoh(2)$~\cite[p.39]{Shimcomp}.
The irrational values of~$x$ that are roots of the irreducible
polynomials of degree~$3$ and~$6$ yield irrational \hbox{$t$-values}
of the same degree.  Hence the \hbox{CM-11} point, being rational,
must have $x = -10099/64$, for which we calculate
the values of~$y$ and~$t$\/ given in the statement of the Corollary.
\Qed

The CM-8 point also has an element of norm~$3$ in its
endomorphism ring, namely $1+\sqrt{-2}$.  The fact that
the CM~points with $x = -10099/64$ have discriminant~$-11$
also follows from the fact that their \hbox{$y$-coordinates}
are conjugates over $\Q(\sqrt{-11})$.
The irrational values of~$x$ have degrees 1 and~$2$ over~$K$\/;
in particular, the cubic irrationalities yield
four \hbox{$\Gal(K/\Q)$-orbits} of CM points.
Taking $c = 2 \cos(2\pi/7)$,
we find that these CM points have \hbox{$x$-coordinates}
$189 c - 19$,
$-(189 c^2 + 567 c + 397)$,
$-2(756 c^2 + 1701 c + 671)$,
$-2(3591 c^2 + 8127 c + 2939)$,
and their $\Gal(K/\Q)$ conjugates.
We determine their CM~fields as the fields of definition
of the points' \hbox{$y$-coordinates}: they are the fields
obtained by adjoining to~$K$\/ the square roots of 
the totally negative algebraic integers
$c^2 + 2c - 7$, $c-6$, $-(3c^2+2c+3)$, and $c^2-2c-11$,
of norms $-167$, $-239$, $-251$, and $-491$ respectively.

\section{The curve $\XXun(3)$, and its Jacobian $\JJun(3)$}

The tables \cite{BK,Cr} of elliptic curves of low conductor
indicate that the curve \hbox{147-B2(J)}, which we identified
with $\XXoh(3)$, is \hbox{$13$-isogenous} with another
elliptic curve over~$\Q$, namely \hbox{147-B1(I)}: $[0,1,1,-114,473]$.
Now \hbox{$13$-isogenies} between elliptic curves over~$\Q$,
though not hard to find (via a rational parametrization of
the classical modular curve $\X_0(13)$), are rare:
this \hbox{$13$-isogeny}, and its twist by $\Q(\sqrt{-7}\,)$
[curves \hbox{147-C1(A)} and \hbox{147-C2(B)}],
are the only examples up to conductor $200$ \cite{BK};\footnote{
  Table 4 in \cite{BK}, which gives all curves of conductor $2^a 3^b$,
  gives several other examples of conductor $20736 = 2^8 3^4$,
  all related by quadratic twists.
  }
even up to conductor~$1000$, the only other examples
are the twists of these \hbox{$13$-isogenies} by $\Q(\sqrt{-3}\,)$,
which appear at conductor~$441$.  The fact that $\XXoh(3)$ admits
a rational \hbox{$13$-isogeny} thus seemed a remarkable coincidence.
Mark Watkins \cite{Watkins} observes that
this curve \hbox{147-B1(I)} has an even more unusual property:
not only is there a \hbox{$13$-isogeny} from this curve to $\XXoh(3)$,
but the kernel of the isogeny consists of points rational over~$K$.
Whereas the classical modular curve $\X_0(13)$, which parametrizes
(generalized) elliptic curves with a rational \hbox{$13$-isogeny},
has genus~$0$, the curve $\X_1(13)$, which parametrizes
(generalized)  elliptic curves with a rational
\hbox{$13$-torsion} point, has genus~$2$.
Thus by Mordell-Faltings $\X_1(13)$ has only finitely many
\hbox{$K$-rational} points.  Hence there are only finitely many
isomorphism classes of elliptic curves defined over~$K$,
let alone over~$\Q$, with a \hbox{$K$-rational}
\hbox{$13$-torsion} point --- and we have found one of them\footnote{
  Not, however, the unique one.  The relevant twist of $\X_1(13)$ by~$K$\/
  is isomorphic with the curve $y^2=4x^6+8x^5+37x^4+74x^3+57x^2+16x+4$,
  which inherits from $\X_1(13)$ a \hbox{$3$-cycle} generated by
  $(x,y) \mapsto (-x^{-1}-1,y/x^3)$.  Rational points related
  by this \hbox{$3$-cycle} and/or the hyperelliptic involution
  parametrize the same curve, with a different choice of generator
  of the torsion subgroup.  The orbit of \hbox{$\Q$-rational} points
  with $x \in {0,-1,\infty}$ yields our curve~$\JJun(3)$.
  A computer search finds that there is at least one other orbit,
  represented by $(x,y)=(17/16,31585/2048)$, which leads to a second
  elliptic curve over~$\Q$ with a \hbox{$K$-rational}
  \hbox{$13$-torsion} point.
  We calculate that this is the curve with coefficients
  $$
  [0, 1, 1, -69819305996260175838254, 7100867435908988429025874812520367]
  $$
  and conductor $8480886141 = 3 \cdot 7^2 \cdot 13 \cdot 251 \cdot 17681$.
  }
by computing the Jacobian of the Shimura modular curve $\XXun(3)$!

These observations are explained by considering $\XXun(3)$.
This curve, like $\XXoh(3)$, is defined over~$\Q$, and
the cyclic cover $\XXun(3) \ra \XXoh(3)$ has degree $(3^3-1)/2 = 13$.
This cover must be unramified, because the only elliptic point
on $\XXoh(3)$ is of order~$3$, which is coprime to~$13$.
Hence $\XXun(3)$ has genus~$1$.  It does not quite follow that
$\XXun(3)$ is isomorphic over~$\Q$ with
the \hbox{$13$-isogenous} elliptic curve \hbox{147-B1(I)},
because $\XXun(3)$ need not have any \hbox{$\Q$-rational} point.
However, the Jacobian $\JJun(3)$ of~$\XXun(3)$ is an elliptic curve,
and is also \hbox{$13$-isogenous} with~$\XXoh(3)$
because the cover \hbox{$\XXun(3) \ra \XXoh(3)$}
induces a map of the same degree from~$\JJun(3)$ to~$\XXoh(3)$.
Since the elliptic curve \hbox{147-B1(I)} is the only one
\hbox{$13$-isogenous} with $\XXoh(3)$ over~$\Q$,
we conclude that it is isomorphic with $\JJun(3)$.
Furthermore, the geometric Galois group of the cover $\XXun(3)/\XXoh(3)$
is canonically $k^*/\{\pm 1\}$,
where $k$ is the residue field of the prime above~$3$ in~$K$.
Working over~$K$, we can choose a generator for this group,
which acts on $\XXun(3)$ by translation by some element of $\JJun(3)$,
and this element is a \hbox{$13$-torsion} point that generates
the kernel of the isogeny from~$\JJun(3)$ to~$\XXoh(3)$.

This is all quite reminiscent of the situation for the classical
modular curves $\X_0(11)$, $\X_1(11)$,
which are \hbox{$5$-isogenous} elliptic curves,
with the kernel of the isogeny $\X_1(11) \ra \X_0(11)$
canonically isomorphic with $\F_{11}^*/\{\pm1\}$
and thus consisting of \hbox{$\Q$-rational} \hbox{$5$-torsion} points.
Also, $\X_1(11)$ has considerably smaller height
(visible in its smaller coefficients) and discriminant than $\X_0(11)$,
as does $\JJun(3)$ compared with $\XXoh(3)$.
The comparison becomes even clearer if we work
with models of these curves minimal over~$K$,
where the additive reduction at the rational prime~$7$
becomes good reduction at the prime $(2-c)$ of~$K$\/ above~$7$.
%
%
%
%

Now for the classical modular curves
the fact that $\X_1(11)$ has a simpler equation than $\X_0(11)$
illustrates a general phenomenon noted in~\cite{Stevens}:
the minimal height in an isogeny class of elliptic curves over~$\Q$
is conjecturally attained by the optimal quotient
of the Jacobian of $\X_1(N)$, not $\X_0(N)$
(unless the \hbox{$\X_1$- and} \hbox{$\X_0$-optimal} quotients coincide).
Does this always happen also for Shimura curves?
We observed the same behavior in several other cases,
one of which appears in the next section.
But there are no extensive tables of optimal quotients
of Jacobians of Shimura curves on which one might test
such a conjecture.  Although Vatsal \cite[Thm.~1.11]{Vatsal}
has proven Stevens' conjecture for curves with a rational
\hbox{$\ell$-isogeny} for prime $\ell\geq 7$,
his methods cannot apply in our setting
(even though all our curves have a suitable isogeny)
because they rely on congruences between \hbox{$q$-expansions}
of modular forms, a tool that is not available to us
in the Shimura-curve setting.

\section{Some other Shimura curves of genus~$1$}

Watkins notes several other examples of curves in~\cite{Cr}
that behave similarly for other choices of cyclic cubic fields~$K$,
and asks whether they, too, can be explained as Shimura modular curves
or Jacobians.  We checked that this is the case for at least one
of these examples, at conductor~$162 = 2\cdot 9^2$.
In this section, we outline this computation and describe
a \hbox{$7$-isogeny} in conductor~$338 = 2 \cdot 13^2$
that should also involve Shimura modular curves.

For the curves of conductor~$162$, we start with a quaternion
algebra over $K_9 = \Q(\cos 2\pi/9)$ ramified at two of the three
real places of $K_9$ and at none of its finite primes.
The resulting modular group $\Gamma(1)$ is again a triangle group,
this time with signature $(2,3,9)$ rather than $(2,3,7)$
(see \cite{Takeuchi}, class~XI).  Since the rational prime~$2$
is inert in~$K_9$, we have modular curves $\XXoh(2)$ and $\XXun(2)$
with geometric Galois group $\PSL_2(\F_8)$ over the rational curve $\Xone$.
We calculate that here $\XXoh(2)$ already has genus~$1$.
Since the Belyi map \hbox{$\XXoh(2) \ra \Xone$} has degree
as low as $8+1=9$, we can find its coefficients with little difficulty
even without resorting to the methods we used to obtain
the equations in Theorem~1.
We place the elliptic points of $\Xone$ at $t=1,0,\infty$.
Then $t$ is a rational function of degree~$9$ on $\XXoh(2)$
with a ninth-order pole, which we use as the base point~$P_0$
to give $\XXoh(2)$ the structure of an elliptic curve.
The zero divisor of~$t$ is then $3D_3$
for some divisor~$D_3$ of degree~$3$.
Hence $D_3-3(P_0)$ is a \hbox{$3$-torsion} divisor on~$\XXoh(2)$.
By computing the divisor of $dt$\/ as in the proof of Proposition~1(i),
we see that the corresponding \hbox{$3$-torsion} point on~$\XXoh(2)$
is also the simple zero of~$t-1$.
In particular, this \hbox{$3$-torsion} point cannot be trivial.
(This fact could also be deduced by noting that
if the \hbox{$3$-torsion} point were trivial then $t$
would be the cube of a rational function on~$\XXoh(2)$,
and so could not have Galois group $\PSL_2(\F_8)$.)
The general elliptic curve with a rational \hbox{$3$-torsion} point
has the form $y^2 + a_1 x y + a_3 y = x^3$,
with the torsion point at $(x,y)=(0,0)$.
Solving for $(a_3\0:a_1^3)$ and the coefficients of~$t$,
 we soon find a model of $\XXoh(2)$ with $(a_1,a_3) = (15,128)$
and $t = (y-x^2-17x)^3 / (2^{14} y)$, and with the double roots of~$t-1$
occurring at the zeros of $(x+9)y + 6x^2 + 71x$ other than $(0,0)$.
Reducing $\XXoh(2)$ to standard minimal form, we find the curve
\hbox{162-B3(I)}: $[1,-1,1,-95,-697]$.  As expected, this curve
attains its \hbox{$3$-adic} potential good reduction over~$K_9$.
The involution $w_2$ is the unique involution of this elliptic curve
that fixes the simple zero $(0,-128)$ of $t-1$.  For instance,
$w_2$ sends the point at infinity to the other \hbox{$3$-torsion} point
$(0,0)$, at which $t=-17^3/2^7$.
Hence $t = -17^3/2^7$ is a CM point on $\Xone$, the unique point
\hbox{$2$-isogenous} with the elliptic point $t=\infty$.
Solving $t(P)=t(w_2(P))$ yields the further CM point
$t = 17^3 5^3 / 2^6$, which must have CM field $K_9(\sqrt{-7}\,)$
because $P$\/ and~$w_2(P)$ are conjugate over $\Q(\sqrt{-7}\,)$.

Besides the known \hbox{$3$-torsion} point,
we find that $\XXoh(2)$ also has a \hbox{$7$-isogeny}
with the curve \hbox{162-B1(G)}: $[1,-1,1,-5,5]$.
It is known (see for instance the ``remarks on isogenies'' in~\cite{BK})
that the isogeny class of this curve is the unique one,
up to quadratic twist, with both a $3$- and a \hbox{$7$-isogeny}
over~$\Q$.
%
%
%
We have already accounted for the \hbox{$3$-torsion} point
using the ramification behavior of the map $\XXoh(2) \ra \Xone$.
The \hbox{$7$-isogeny} is again explained using $\XXun$.
This time the cyclic cover $\XXun(2) \ra \XXoh(2)$
has degree $2^3 - 1 = 7$, and the only elliptic point of $\XXoh(2)$
is the simple zero of~$t-1$, which has order~$2$.
Since $\gcd(2,7)=1$, the cyclic cover is again unramified,
and $\JJun(2)$ must be an elliptic curve
\hbox{$7$-isogenous} with~$\XXoh(2)$.  Hence $\JJun(2)$
is isomorphic with the elliptic curve \hbox{162-B1(G)}: $[1,-1,1,-5,5]$.
Also as before, the kernel of the isogeny must consist of points 
rational over~$K_9$, and again Watkins confirms~\cite{Watkins}
that this is the case for this curve.
This time it turns out $\JJun(2)$ is the unique
elliptic curve over~$\Q$ with the correct torsion structure:
a \hbox{$3$-torsion} point over~$\Q$, and a \hbox{$7$-torsion} point
over~$K_9$ that is proportional to its Galois conjugates.
Again we note that it is $\JJun(2)$, not $\XXoh(2)$,
that has the smaller discriminant and height.

The next case is the cyclic cubic field
$K_{13}$ of discriminant $169=13^2$.
%
Once more we use a quaternion algebra over this field
that is ramified at two of its three real places
and at none of the finite primes.  By Shimizu's formula
(\cite[Appendix]{Shimizu}, quoted in \cite[p.207]{Takeuchi}),
the resulting curve $\Xone$ has hyperbolic area
$$
\frac1{16\pi^6} d_{K_{13}}^{3/2}\zeta_{K_{13}}(2)
= -\frac12 \zeta_{K_{13}}(-1).
$$
Since $K_{13}$ is cyclotomic, we can compute
$\zeta_{K_{13}}(2)$ or $\zeta_{K_{13}}(-1)$
by factoring $\zeta_K$ as a product of the Riemann zeta function
and two Dirichlet \hbox{$L$-series}.
We find that $\Xone$ has normalized hyperbolic area~$1/6$.
Since $\Gamma(1)$ is not on Takeuchi's list
of triangle groups~\cite{Takeuchi}, the curve $\Xone$
must have at least four elliptic points,
or positive genus and at least one elliptic point.
The only such configuration that attains an area as small as $1/6$
is genus~zero, three elliptic points of order~$2$, and one of order~$3$.
Again we use the prime of~$K$\/ above the inert rational prime~$2$
to construct modular curves $\XXoh(2)$ and $\XXun(2)$ with
geometric Galois group $\PSL_2(\F_8)$ over the rational curve $\Xone$.
In the \hbox{degree-$9$} cover $\XXoh(2) \ra \Xone$,
the elliptic point of order~$3$ has $3$ triple preimages,
and each of the \hbox{order-$2$} points
has one simple and four double preimages.
Hence $\XXoh(2)$ has genus~$1$ by Riemann-Hurwitz.
Again the cyclic cover $\XXun(2) \ra \XXoh(2)$ has degree~$7$,
relatively prime to the orders of the elliptic points on~$\XXoh(2)$,
so $\XXun(2)$ also has genus~$1$.

In this setting it is not clear that either $\XXoh(2)$ or~$\XXun(2)$
is an elliptic curve over~$\Q$, since neither curve is forced
to have a \hbox{$\Q$-rational} divisor of degree~$1$.
(The elliptic points of order~$2$ may be Galois conjugates,
not individually rational over~$\Q$.)
But we can still consider the Jacobians $\JJoh(2)$ and $\JJun(2)$,
which are elliptic curves \hbox{$7$-isogenous} over~$\Q$,
with the kernel of the \hbox{$7$-isogeny} $\JJun(2) \ra \JJoh(2)$
consisting of points rational over~$K_{13}$.  Watkins suggests,
by analogy with the cases of conductor $147$ and~$162$,
that this \hbox{$7$-isogeny} should connect
curves of conductor $2\cdot 13^2 = 338$.
According to the tables of~\cite{Cr}, there are
three pairs of \hbox{$7$-isogenous} curves of conductor~$338$.
Watkins computes that of the $6$ curves involved in these isogenies,
only \hbox{338-B1}: $[1,-1,1,137,2643]$ has \hbox{$7$-torsion} points
rational over~$K_{13}$, and thus proposes this curve as $\JJun(2)$,
and the \hbox{$7$-isogenous} curve \hbox{338-B2}:
$[1, -1, 1, -65773, -6478507]$ as $\JJoh(2)$.
This must be correct, since these Jacobians should again
have multiplicative reduction at~$2$ and good reduction
at all primes of~$K_{13}$ other than~$2$, whence their conductor
must be $2\cdot 13^2$ (there being no curve of conductor~$2$ over~\Q).
But here I have not obtained an explicit rational function
or even determined the cross-ratio of the elliptic points on~$\Xone$,
which would be quite a demanding computation
using the methods of~\cite{Shimcomp}.

\vspace*{2ex}

\noindent
{\bf Acknowledgements.}
This work was supported in part by NSF grants DMS-0200687 and DMS-0501029.
Symbolic and numerical computations were assisted by
the computer packages {\sc pari/gp}~\cite{Pari} and {\sc maxima}.
I am also grateful to Mark Watkins for the reference to~\cite{Vatsal}
and other relevant correspondence~\cite{Watkins},
to Watkins and the referees for careful readings of and corrections to
earlier drafts of this paper,
to the referees for suggesting several other bibliographical sources,
and to Benedict Gross for the remarks in footnote~\ref{footnote:carayol}.

\end{document}